\documentclass[11pt]{amsart}
\usepackage{amsthm,amsfonts,amsmath,amssymb,latexsym,epsfig,mathrsfs,yfonts,marvosym}
\usepackage{graphicx,color}
\usepackage{url}
\usepackage[colorlinks=true]{hyperref}
\usepackage[utf8]{inputenc}
\usepackage[T1]{fontenc}

\newtheorem{remark}{Remark}[section]

\newtheorem{corollary}{Corollary}
\newtheorem{example}{Proposition - Example}

\def\<{\langle}

\def\>{\rangle}

\def\a{\alpha}

\newcommand{\grad}{\mathop{\mathrm{grad}}\nolimits}

\begin{document}

\title[Non trivial examples of coupled equations]{Non trivial examples of coupled equations for K\"ahler metrics and Yang-Mills connections}
\author{Julien Keller and Christina W. T{\o}nnesen-Friedman}
\thanks{This work was partially supported by a grant from the Simon's Foundation (208799 to Christina T{\o}nnesen-Friedman)}
\address{Julien Keller \\ Centre de Math\'ematiques et Informatique, Universit\'e de Marseille, \\ Technop\^ole de Ch\^ateau Gombert, 39 rue F. Joliot-Curie \\ 13453 Marseille Cedex 13 \\ France} \email{jkeller@cmi.univ-mrs.fr} 
\address{Christina W. Tonnesen-Friedman\\ Department of Mathematics\\ Union
College\\ Schenectady\\ New York 12308\\ USA } \email{tonnesec@union.edu}

\begin{abstract}
We provide non trivial examples of solutions to the system of coupled equations introduced by M. Garc\'ia-Fern\'andez for the uniformization problem of a triple $(M,L,E)$ where $E$ is a holomorphic vector bundle over a polarized complex manifold $(M,L)$, generalizing the notions of both constant scalar curvature K\"ahler metric and Hermitian-Einstein metric.
\end{abstract}
\maketitle

\section{Introduction}
In his Ph.D thesis \cite{FG} (see also \cite{AGFGP}), M. Garc\'ia-Fern\'andez introduced a natural system of equations, called ``coupled equations'' that is related to a very simple data,
composed by a K\"ahler manifold, a K\"ahler class and a vector bundle on the underlying manifold. Let us recall briefly how this system appears through a moment map construction. \par 
Let $(M,\omega)$ be a compact symplectic manifold, $G$ a compact Lie group and $\pi:E\rightarrow M$ a smooth principal $G$-bundle on $M$.
Let $\mathcal{J}$ be the space of almost complex structures compatible with $\omega$, and $\mathcal{A}(E)$ be the space
of connections on $E$. One can consider the extended gauge group $\widetilde{\mathcal{G}}$ of $E$. It is the group of automorphisms
of $E$ that cover the hamiltonian symplectomorphisms $\mathcal{H}$ of $M$. We mean that for such an automorphism $a$,  $\pi \circ a= h \circ \pi$ where
$h\in\mathcal{H}$ is a symplectomorphism of $M$, and in particular this group  $\widetilde{\mathcal{G}}$ contains the gauge group of $E$.
It is possible to check that there is a surjection map $p:\widetilde{\mathcal{G}}\mapsto \mathcal{H}$ \cite[Section 2.2]{FG}. Using this surjection, one can define
an action of $a\in \widetilde{\mathcal{G}}$ on $\mathcal{J}\times \mathcal{A}(E)$, by 
$$a \cdot(J,A)=(p(g)_* J,g\cdot A)=(dp(g)\circ J\circ (dp(g))^{-1},g\circ A\circ g^{-1}).$$
Using the symplectic structures $\omega_{\mathcal{J}},\omega_{\mathcal{A}(E)}$ on  $\mathcal{J}$ and $\mathcal{A}(E)$, one can fix a symplectic form on $\mathcal{J}\times \mathcal{A}(E)$ by 
\begin{equation}\label{symplform}
\Omega=\alpha_0\omega_{\mathcal{J}} +\frac{\alpha_1}{(n-1)!}\omega_{\mathcal{A}(E)},
\end{equation}
for $\alpha_0,\alpha_1$ two real constants.
Of course this choice is not canonical a priori, but this is probably the simplest one. In this setting, there is a moment map $\mu : \mathcal{J}\times \mathcal{A}(E) \rightarrow Lie(\widetilde{\mathcal{G}})^*$ associated to the action of  $\widetilde{\mathcal{G}}$ and the symplectic structure $\Omega$  we have just described \cite[Proposition 2.3.1]{FG}. \par
Let us assume from now that $M$ has a K\"ahler structure, $\omega$ is a K\"ahler form,. We can restrict our attention to $\mathcal{J}^i\subset \mathcal{J}$ the 
set of integrable almost-complex structures that are compatible with $\omega$, and $\mathcal{A}^{1,1}_J(E)$, the subset of compatible connections $A$ on $E$ such that $F_A^{0,2}=F^{2,0}_A=0$ with respect to the structure $J$. There is a $\widetilde{\mathcal{G}}$-invariant complex submanifold $S\subset \mathcal{J}\times \mathcal{A}(E)$ consisting of pairs $(J,A)$ such that $ J\in \mathcal{J}^i$ and $A\in \mathcal{A}^{1,1}_J(E)$.
 There is a K\"ahler form on the non singular part of this complex submanifold $S$
obtained by restriction of $\Omega$ and thus an induced moment map $\mu$ from the $\widetilde{\mathcal{G}}$-action by holomorphic isometries. \\
A pair  $(J,A)\in S$ is solution of the following system of two equations, namely the coupled equations
\begin{equation}\label{cpled}
{\Big\{}\begin{array}{ccc}
    \Lambda_\omega F_A &=& z,  \\
  \alpha_0 Scal(g_J)+\alpha_1\Lambda^2_{\omega}(F_A\wedge F_A)&=&\alpha_2.
  \end{array}
\end{equation}
if its orbit belongs to the K\"ahler reduction $\mu^{-1}(0)/\widetilde{\mathcal{G}}$.
Here $g_J$ is the metric $g_J=\omega(\cdot, J\cdot)$ induced by the K\"ahler 
form $\omega$ and the structure $J\in \mathcal{J}^i$, $Scal(g_J)$ the scalar 
curvature of $g_J$, $\Lambda_\omega$ the contraction operator with respect to 
the form $\omega$, $F_A$ the curvature of $A$, and $\alpha_2$ a constant 
dependent on $\alpha_0$, $\alpha_1$, the K\"ahler class, and topological constants. Likewise
$z$ is determined by the K\"ahler class, and topological constants (see Remark 1.2 in \cite{AGFGP}).

Of course when one considers a trivial bundle $E$, it turns out that the coupled system can be solved by finding a constant scalar curvature K\"ahler metric and a flat connection. From Fujiki and Donaldson's work, it is 
well-known that the constant scalar curvature K\"ahler equation (cscK equation 
in short) appears as prescribing a zero of the moment map induced by the action 
of the hamiltonian symplectomorphisms $\mathcal{H}$ on the integrable complex 
structures $\mathcal{J}^i$, see \cite{D4} for details. On the other hand, the 
first equation of (\ref{cpled}) appears naturally when one is considering the 
action of the gauge group of $E$ \cite[Chapter 6]{DK} in view of the 
Kobayashi-Hitchin correspondence. Certainly, the motivation to study the coupled
equations (\ref{cpled}) is coming from the natural question: what is a good 
moduli problem for a tuple $(M,L,E)$ where $(M,L)$ is a compact polarised manifold with a K\"ahler class $2\pi c_1(L)$ and $E$ is a holomorphic $G^\mathbb{C}$-bundle over $M$ ? It is also natural to wonder what are the natural geometric perturbations of the cscK equation, and to analyze its perturbations in terms of K-stability \cite[Chapter 4]{FG}, in order to test the so-called Yau-Tian-Donaldson conjecture. 

\medskip

In order to do so, one is lead to construct non trivial examples 
($\alpha_0\alpha_1\neq 0$) of solutions to the coupled equations (\ref{cpled}) 
for the tuple $(M,L,E)$ chosen as before. M. Garc\'ia-Fern\'andez showed that it
is possible to obtain examples by deformations (when $\frac{\alpha_1}{\alpha_0}$
is small enough)  of a manifold $M$ that carries a constant scalar curvature
K\"ahler metric and a Hermitian-Yang-Mills holomorphic vector bundle, if the 
automorphism group of $M$ is finite \cite[Theorem 3.2.4]{FG}. The proof reduces 
to an implicit function theorem and the assumption on the automorphism group
allows to invert the linearization operator. Also, when $M$ has complex
dimension 1, the term $F_A\wedge F_A$ vanishes and one can provide solutions 
to the coupled equations by considering the Kobayashi-Hitchin correspondence 
for holomorphic bundles, i.e Mumford polystable bundles on the complex curve. 
On compact homogeneous K\"ahler-Einstein surfaces, examples can be provided by 
considering anti-self-dual connections.  In higher dimensions, examples can be 
found by using projectively flat bundles over a manifold with constant scalar 
curvature metric and satisfying a natural topological condition. As one can 
remark, all these examples are very specific since they hold on manifolds that 
carry a constant scalar curvature K\"ahler metric. It is natural to wonder if 
one can find new examples of solutions to the coupled equations on 
complex manifolds 
such that there is no solution to the cscK equation in the class $2\pi c_1(L)$. The main goal of this paper is to construct such examples over a ruled surface and a ruled threefold.

\medskip
{ \small \noindent We would like to thank Mario Garc\'ia-Fern\'andez for very useful discussions. 
We are also thankful for numerous helpful conversations with Vestislav Apostolov, David Calderbank, and Paul Gauduchon.
C. T{\o}nnesen-Friedman would like to thank them
for involving her in their project on admissible manfolds long time ago. The results of that project play a crucial technical role in this paper.}

 \section{Examples of solutions to coupled equations on Hirzebruch type ruled 
surfaces}
Let us consider a ruled manifold of the form
$M = {\mathbb P}({\mathcal O} \oplus {\mathcal L}) \rightarrow \Sigma$,
where $\Sigma$ is a compact Riemann surface, ${\mathcal L}$ is a  
holomorphic line bundle of degree $k\in {\mathbb Z}^{*}_+$ on $\Sigma$, and 
${\mathcal O}$ is the trivial holomorphic line bundle.
Let $g_{\Sigma}$ be the
K\"ahler metric
on $\Sigma$ of constant scalar curvature $2s_{\Sigma}$, with K\"ahler form
$\omega_{\Sigma}$, such that
$c_{1}({\mathcal L}) = [\frac{\omega_{\Sigma}}{2 \pi}]$.
Let ${\mathcal K}_\Sigma$ denote the canonical bundle of $\Sigma$. Since $c_1({\mathcal K}_\Sigma^{-1}) = [\rho_\Sigma/2\pi]$, where $\rho_\Sigma$ denotes the Ricci form, we have the relation $s_{\Sigma}= 2(1-h)/k$, where $h$
denotes the genus of $\Sigma$.

The natural $\mathbb{C}^*$-action on ${\mathcal L}$ 
extends to a holomorphic
$\mathbb{C}^*$-action on $M$. The open and dense set $M_0$ of stable points with 
respect to the
latter action has the structure of a principal $\mathbb{C}^*$-bundle over the stable 
quotient.
The hermitian norm on the fibers induces via a Legendre transform a function
$\mathfrak{z}:M_0\rightarrow (-1,1)$ whose extension to $M$ consists of the critical manifolds
$E_{0}:=\mathfrak{z}^{-1}(1)=P({\mathcal O} \oplus 0)$ and 
$E_{\infty}:= \mathfrak{z}^{-1}(-1)=P(0 \oplus {\mathcal L})$.

These zero and infinity sections, $E_{0}$ and $E_{\infty}$, of $M \rightarrow 
\Sigma$ have the property that $E_{0}^{2} = k$ and $E_{\infty}^{2} = -k$, 
respectively. If $C$ denotes a fiber of the ruling $M \rightarrow 
\Sigma$, then $C^{2}=0$, while $C \cdot E_{i} =1$ for both, $i=0$ and 
$i=\infty$. Any real cohomology class in the two dimensional space 
$H^{2}(M, {\mathbb R})$ may be written as a linear combination of 
(the Poincar\'e 
duals of) $E_{0}$ and $C$, 
$$
m_{1} E_{0} +m_{2} C\, .
$$
Thus, we may think of $H^{2}(M, {\mathbb R})$ as ${\mathbb R}^2$, 
with coordinates $(m_{1},m_{2})$. The K\"ahler cone ${\mathcal K}$ may 
be identified with ${\mathbb R}_{+}^2=\{ (m_{1}, m_{2})\,  | \; m_{1} > 0, 
m_{2} > 0 \}$ (see \cite{fujiki} or Lemma 1 in \cite{t-f}).
To calculate $m_{1}$ and $m_{2}$ for a real cohomology class $\Gamma 
\in H^{2}(M, {\mathbb R})$ 
it is useful to notice that we have
\begin{eqnarray}\label{inter}
\Big{\{}\begin{array}{l}
\int_{E_{0}} \Gamma = \Gamma \cdot E_{0} = (m_{1} E_{0} +m_{2} C) \cdot E_{0}   =   k m_{1} + m_{2},\\
\int_{C} \Gamma = \Gamma \cdot C =(m_{1} E_{0} +m_{2} C) \cdot C  =  m_{1}.
\end{array}
\end{eqnarray}
Thus, we get $m_{1} = \int_{C} \Gamma$ and $m_{2}= \int_{E_{0}} \Gamma - 
k \int_{C} \Gamma$.

\medskip We shall use the techniques developed in \cite{acg,acgt} to build a solution to the coupled equations (\ref{cpled}) and check that the manifold $M$ does not carry a constant scalar curvature K\"ahler metric. 
To build the so-called admissible metrics \cite{acgt} on $M$ we 
proceed as follows. Let  $\theta$ be a connection one form for the 
Hermitian metric on $M_0$, with curvature
$d\theta = \omega_\Sigma$. Let $\Theta$ be a smooth real function with 
domain containing
$(-1,1)$. Let $x$ be a real number such that $0 < x < 1$.
Then an admissible K\"ahler metric
is given on $M_0$ by
\begin{equation}\label{metric}
g  =  \frac{1+x \mathfrak{z}}{x} g_\Sigma
+\frac {d\mathfrak{z}^2}
{\Theta (\mathfrak{z})}+\Theta (\mathfrak{z})\theta^2\,
\end{equation}
with K\"ahler form 
\begin{equation}
\omega =  \frac{1+x \mathfrak{z}}{x}\omega_\Sigma
+d\mathfrak{z}\wedge \theta\,. \label{kf}
\end{equation}
The complex structure yielding this
K\"ahler structure is given by the pullback of the base complex structure
along with the requirement 
\begin{equation}\label{complex}
Jd\mathfrak{z} = \Theta \theta  
\end{equation}
 The function $\mathfrak{z}$ is 
hamiltonian
with $K= J\grad \mathfrak{z}$ a Killing vector field. Observe that $K$ 
generates the circle action which induces the holomorphic
$\mathbb{C}^*$- action on $M$ as introduced above.
In fact, $\mathfrak{z}$ is the moment 
map on $M$ for the circle action, decomposing $M$ into 
the free orbits $M_{0} = \mathfrak{z}^{-1}((-1,1))$ and the special orbits 
$\mathfrak{z}^{-1}(\pm 1)$. Finally, $\theta$ satisfies
$\theta(K)=1$.
In order that $g$ (be a genuine metric and) extend to all of $M$,
$\Theta$ must satisfy the positivity and boundary
conditions
\begin{align}
\label{positivity}
(i)\ \Theta(\mathfrak{z}) > 0, \quad -1 < \mathfrak{z} <1,\quad
(ii)\ \Theta(\pm 1) = 0,\quad
(iii)\ \Theta'(\pm 1) = \mp 2.
\end{align}
The last two of these are together necessary and sufficient for
the compactification of $g$.
Define a function $F(\mathfrak{z})$ by the formula 
\begin{equation}\label{theta}
\Theta(\mathfrak{z})= \frac{F(\mathfrak{z})}{(1+x 
\mathfrak{z})} 
\end{equation}
Since $(1+x 
\mathfrak{z})$ is positive for $-1<\mathfrak{z}<1$, conditions
\eqref{positivity}
imply the following equivalent conditions on $F(\mathfrak{z})$:
\begin{align}
\label{positivityF}
(i)\ F(\mathfrak{z}) > 0, \quad -1 < \mathfrak{z} <1,\quad
(ii)\ F(\pm 1) = 0,\quad
(iii)\ F'(\pm 1) = \mp 2(1 \pm x).
\end{align}

The volume form of $g$ in \eqref{metric} is given by
\begin{equation*}
d\mu_{g} = \frac{\omega \wedge \omega}{2} = \frac{1+x\mathfrak{z}}{x} 
\omega_{\Sigma} \wedge d\mathfrak{z} \wedge \theta \,,
\end{equation*}
while the Ricci form is given by
$$
\rho_{g} = \left( s_{\Sigma} - \frac{F'(\mathfrak{z})}{2(1+ x \mathfrak{z})} \right) 
\omega_{\Sigma} - 
d\left( \frac{F'(\mathfrak{z})}{2(1+ x \mathfrak{z})} \right) 
\wedge \theta \, ,$$
and the scalar curvature is given by
$$ Scal(g) = \frac{2s_{\Sigma} x}{1+x\mathfrak{z}} - \frac{F''(\mathfrak{z})}{1+x\mathfrak{z}}\, .$$
The calculations of these geometrical terms can be found in 
\cite{acg}.

Now, since $E_{0} = \mathfrak{z}^{-1}(1)$ and $k = c_{1}({\mathcal L}) = 
[\frac{\omega_{\Sigma}}{2 
\pi}]$ we have that 
$$\int_{E_{0}}[\omega] = \frac{1+x}{x} 
\int_{\Sigma}\omega_\Sigma = \frac{2\pi k (1+x)}{x}.$$ It is also easy 
to see that 
$$\int_{C}[\omega] = \int_{0}^{2\pi}\int_{-1}^{1}d\mathfrak{z}\wedge 
dt = 4 \pi,$$ where $t \in [0,2\pi]$ is a fibre coordinate of 
$M_{0} \rightarrow \Sigma \times (-1,1)$ in a gauge chosen such that 
the connection form $\theta$ has no $d\mathfrak{z}$ components. Therefore, from (\ref{inter}), we have that 
\begin{equation}\label{class1}
[\omega] = 4\pi E_{0}+ \frac{2\pi(1-x)k}{x} C \, .
\end{equation}
and we fix $x=\frac{k}{k+k'}$ with $k'\in \mathbb{Z}^*_+$.
\medskip 
The $(1,1)$ form 
\begin{align*}
 \rho_{g} - \frac{Scal(g)}{4} \omega = \left(\frac{2 s_{\Sigma}x(1+x \mathfrak{z}) + 
F''(\mathfrak{z}) (1+ x \mathfrak{z}) - 2 x F'(\mathfrak{z})}{4(1+x \mathfrak{z})^{2}}\right) \\
   &\hspace{-3cm}\times \left(\frac{1+x \mathfrak{z}}{x}\omega_\Sigma
- d\mathfrak{z}\wedge \theta\right)
\end{align*}
is traceless and therefore 
anti-self-dual. Using this, we easily check that the form
$$
\a := \frac{x^{2}}{(1+x \mathfrak{z})^{2}}\left( \frac{1+x \mathfrak{z}}{x}\omega_\Sigma
- d\mathfrak{z}\wedge \theta \right),
$$
which does not depend on $F(\mathfrak{z})$, is both closed and anti-self-dual.
Since the second Betti number $b_{2}(M)$ of our ruled surface $M$ is two, 
while the signature, $\sigma$, is zero, a basis for the vector space of 
harmonic real $(1,1)$-forms on $(M,g)$ would be given by 
$\{ \omega, \a\}$.
Now, 
$$\int_{E_{0}}[\a] = \frac{x}{(1+x)}\int_{\Sigma}\omega_\Sigma=
\frac{2\pi k x}{(1+x )}$$
and 
$$\int_{C}[\a] =
-x^{2}\int_{0}^{2\pi}\int_{-1}^{1}(1+x \mathfrak{z})^{-2}d\mathfrak{z}\wedge 
dt = 
\frac{-4\pi x^{2}}{1-x^{2}}.$$ 
We therefore have that
$$
[\a] =\frac{2\pi x}{1-x^2} \left( -2x E_0 +k(1+x)C\right) \, .
$$
Consider the $(1,1)$ form $\gamma_{a,b} = a \, \omega + b\, \a$ for some constants  $a,b$.
Since 
\begin{equation}\label{class2}
\left[\frac{\gamma_{a,b}}{2\pi}\right] = \frac{2(a(1-x^{2}) - b 
x^{2})}{1-x^{2}}\, E_{0} + \frac{k(a(1-x)^{2} + bx^{2})}{x(1-x)}\, C
\end{equation}
it is easy to see that for appropriate choices of $a$ and $b$,
$[\frac{\gamma_{a,b}}{2\pi}]$ is an integral class and thus 
$\gamma_{a,b}$ may be viewed as the curvature form $F_{A}$ of some 
connection $A$ on some stable vector, $E$, bundle over $M$ (of course in that case $E$ is a line 
bundle). Actually any choice of $a\in \mathbb{Z}$ and $b$ integer multiple of 
$\frac{(2k+k')k'}{k^2}$  will imply that  $[\frac{\gamma_{a,b}}{2\pi}]$ belongs to $H^{2}(M, {\mathbb Z})$.
\par We easily calculate that \begin{equation} \Lambda_{\omega} \gamma_{a,b} = 2a, \label{1cpled} \end{equation} which corresponds
to the first equation of  the system (\ref{cpled}). 
The second equation of the coupled equations (\ref{cpled}) in the variables $(g,\omega, J)$ on $M$ and 
$\gamma_{a,b}$ on $E$, corresponds to
\begin{equation}\label{coupled}
\textstyle \alpha_{0} Scal(g) + \alpha_{1}  \Lambda_\omega^{2}\left(\gamma_{a,b}\wedge 
\gamma_{a,b}\right) = \alpha_{2},
\end{equation}
for some constants $\alpha_{0},\alpha_{1}, \alpha_{2} \in {\mathbb R}$.
It is straightforward to verify that \begin{equation}\label{lambda2}
 \gamma_{a,b}\wedge \gamma_{a,b} = 
2\left(a^{2} - 
\frac{b^{2}x^{4}}{(1+x\mathfrak{z})^{4}}\right) \, d\mu_{g}\, ,
 \end{equation} 
and since
$\Lambda^{2}_\omega d\mu_{g} = 2$ (see e.g. 2.77 in \cite{besse}), we have that 
$$\textstyle \Lambda^{2}_\omega \left(\gamma_{a,b}\wedge \gamma_{a,b}\right) = 
4\left(a^{2} - 
\frac{b^{2}x^{4}}{(1+x\mathfrak{z})^{4}}\right).$$ Assuming that $(g,\omega,J)$ is 
admissible, and hence determined by $F(\mathfrak{z})$ satisfying 
\eqref{positivityF}, we get that \eqref{coupled} is equivalent to
$$ \alpha_{0}\left(\frac{2s_{\Sigma} x}{1+x\mathfrak{z}} - \frac{F''(\mathfrak{z})}{1+x\mathfrak{z}}\right) + 
4\alpha_{1}\left(a^{2} - 
\frac{b^{2}x^{4}}{(1+x\mathfrak{z})^{4}}\right) = \alpha_{2}.$$
Unless $b=0$ we have that $\alpha_{0}$ must 
be non-zero.  Otherwise, if $b=0$, we only get a trivial solution since $M$ admits no 
constant scalar curvature K\"ahler metrics.  We therefore arive at the following ODE
\begin{equation}\label{ode}
F''(\mathfrak{z}) = 2 s_{\Sigma} x + \frac{4 \alpha_{1}}{\alpha_{0}}\left(a^{2} - 
\frac{b^{2}x^{4}}{(1+x\mathfrak{z})^{4}}\right)(1+x \mathfrak{z}) - \frac{\alpha_{2}}{\alpha_{0}}(1+x \mathfrak{z}).
\end{equation} 
Integrating twice we see that this has a solution, satisfying 
\eqref{positivityF}, if and only if
\begin{equation}\label{ratio}   
\frac{\alpha_{1}}{\alpha_{0}} = 
\frac{-(1-x^{2})^{2}(2-s_{\Sigma}x)}{8b^{2}x^{4}}\end{equation}
and
\begin{equation}\label{ratio2} \frac{\alpha_{2}}{\alpha_{0}}= \frac{3 b^{2}x^{4}(2+s_{\Sigma}x) - 
a^{2}(2-s_{\Sigma}x)(1-x^{2})^{2}}{2 b^{2} x^{4}}. 
\end{equation}
In that case we find a unique solution
\begin{equation}\label{F}
F(\mathfrak{z}) = \frac{(1-\mathfrak{z}^{2}) (x^{2}(2+s_{\Sigma}x)\mathfrak{z}^{2} + 8x\mathfrak{z} + 4 + 
2x^{2} - s_{\Sigma} x^{3})}{4(1+x \mathfrak{z})}
\end{equation}
and then
\begin{equation} \label{scal}
Scal(g)=\frac{3 (2 + s_{\Sigma} x)}{2}  - \frac{(2 - s_{\Sigma} x) 
(1 - x^2)^2}{2 (1 + x \mathfrak{z})^4}\,\end{equation} Since $s_{\Sigma} x < 2$, $Scal(g)$ 
is not an affine function of $\mathfrak{z}$ and hence not
an extremal K\"ahler metric \cite{acg} (remark that the bundle $\mathcal{O}\oplus \mathcal{L}\rightarrow \Sigma$ is not Mumford polystable so $M$ does not have a constant scalar curvature K\"ahler metric). For the same reason, we also notice that
$\frac{\alpha_{1}}{\alpha_{0}} < 0$ and thus the form $\Omega$ defined by (\ref{symplform}) is symplectic and not K\"ahler.
Finally, with notations above,  we set
\begin{equation} \label{fix}
x=\frac{k}{k+k'}, \;\; a=k_1, \;\; b=\frac{k_2(2k+k')k'}{k^2},
\end{equation}
in the previous equations with $k'\in \mathbb{Z}^*_+$, $k_1\in \mathbb{Z}$, 
$k_2\in \mathbb{Z}^*$. Furthermore, the condition $F(\mathfrak{z})>0$ is satisfied 
for all $-1<\mathfrak{z}<1$. Actually, if we set 
$f(\mathfrak{z})=x^2(2+s_\Sigma x)\mathfrak{z}^2+8x\mathfrak{z}+4+2x^2-s_\Sigma x^3$, then $f(-1)= 4(1-x)^2>0$
and  $f(1)=4(x+1)^2>0$.
If $(2+s_\Sigma x)\leq 0$ it is then straightforward to see that 
$f(\mathfrak{z})>0$ for $-1 \leq \mathfrak{z} \leq 1$. Since $f'(-1) = 2 x (4 - 2 x - 
s_{\Sigma} x^2) >0$ the same can be concluded 
if $(2+s_\Sigma x)> 0$.
Eventually, in both cases, $F$ is strictly positive.\\
Thus, we  have obtained that conditions (\ref{positivityF}), equations  (\ref{coupled}) and (\ref{1cpled}) are all satisfied and the system of coupled equation admits a solution in integral classes from (\ref{class1}) and (\ref{class2}). This leads to the following result. 
\begin{example}\label{thm}
 Assume that $M = {\mathbb P}({\mathcal O} \oplus {\mathcal L}) \rightarrow \Sigma$ is a ruled manifold with  $\mathcal{L}$ of degree $k\in \mathbb{Z}^*_+$. Fix  $k'\in \mathbb{Z}^*_+$, $k_1\in \mathbb{Z}$, $k_2\in \mathbb{Z}^*$. Consider the integral classes  
$L:=2E_0+k'C$ and $E:=2(k_1 -k_2)E_0+(k_1k'+k_2(2k+k'))C$. 
Then,  there exists an admissible  K\"ahler metric $\omega\in 2\pi c_1(L)$, 
a complex structure $J$ defined by (\ref{complex}),(\ref{theta}),(\ref{F}), 
and
a connection $A\in \mathcal{A}_J^{1,1}(E)$
such that the triple $(\omega,J,A)$ is a solution to the coupled equations (\ref{cpled}). 
The constants $(\alpha_0,\alpha_1,\alpha_2)$ satisfy 
$\frac{\alpha_1}{\alpha_0}=-\frac{(2-s_\Sigma)k+2k'}{8k_2^2(k+k')}<0$. Furthermore there is no constant scalar curvature K\"ahler metric in $2\pi c_1(L).$
\end{example}
Note that once the bundle $E$ and the class $L$ are fixed as in our proposition, 
the solution $(\omega,J)$ is unique in the set of 
admissible K\"ahler metrics, up to automorphisms. 
\subsection{About the Calabi-Yang-Mills functional \label{CYMsection}}
We are now going to consider the coupled equations from a variational point of view but in a sligthly different setup than in \cite{FG} and \cite{AGFGP} where the constants $\alpha_i$ are all positive.  First of all, using the fact that for any $A\in A^{1,1}(E)$ \begin{equation}
                     |F_A|^2=|\Lambda_\omega F_A|^2-\frac{1}{2}\Lambda_\omega^2( F_A\wedge F_A) \label{identity1}
                    \end{equation}
 the system of coupled equations is equivalent to 
\begin{equation}\label{cpled2}
{\Big\{}\begin{array}{ccl}
    \Lambda_\omega F_A &=& z,  \\
  \alpha_0 Scal(g_J)-2\alpha_1\vert F_A\vert^2&=&\alpha_2-2\alpha_1 \vert z\vert^2,
  \end{array}
\end{equation}
and we shall consider the case  $\alpha_0\neq 0$. It is natural to introduce the Calabi-Yang-Mills type functional 
$$CYM(g,A)=\int_M \left(Scal(g)-2\frac{\alpha_1}{\alpha_0}|F_A|^2\right)^2 
d\mu_g +  \Vert F_A\Vert^2,$$
and the constant 
$$\alpha_2=\frac{4\pi \alpha_0}{(n-1)!}\frac{\langle c_1(M)\cup 
[\omega]^{n-1},[M]\rangle}{Vol_M([\omega])}  + \frac{2\alpha_1}{(n-2)!}\frac{\langle c(E)\cup [\omega]^{n-2},[M]\rangle}{Vol_M([\omega])}$$
which appears when one is integrating over the manifold $M$ the second equation of the system  (\ref{cpled}). Here we assume that the complex structure is fixed and $g$ varies among K\"ahler metrics with a fixed K\"ahler class $2\pi c_1(L).$

\begin{remark}
 Our functional $CYM$ differs from \cite{FG}. Indeed we choose our functional such that, up to a renormalization, it is globally invariant if we do the change of metric $\omega \rightarrow t\omega$ and change accordingly the constants $\alpha_0,\alpha_1$ by $t\alpha_0,t^2\alpha_1$ in (\ref{cpled}). Furthermore, when $\alpha_1\rightarrow 0$, it reduces to precisely the sum of the Calabi functional and the Yang-Mills functional. \end{remark}

Now, with (\ref{identity1}),  we get that
\begin{eqnarray*} 
CYM(g,A)&=& \int_M \left(Scal(g)-2\frac{\alpha_1}{\alpha_0} |F_A|^2 -\frac{\alpha_2}{\alpha_0}+2\frac{\alpha_1}{\alpha_0}|z|^2\right)^2 d\mu_g \\
&&+ \Vert F_A\Vert^2_{L^2} \\
 && + 2\left(\frac{\alpha_2}{\alpha_0}-2\frac{\alpha_1}{\alpha_0}|z|^2\right)\int_M ( Scal(g)- 2\frac{\alpha_1}{\alpha_0}|F_A|^2)d\mu_g \\
&&- \left(\frac{\alpha_2}{\alpha_0}-2\frac{\alpha_1}{\alpha_0}|z|^2 \right)^2  Vol_M([\omega]),\\
&=& \big\Vert  Scal(g)-2\frac{\alpha_1}{\alpha_0} |F_A|^2- \frac{\alpha_2}{\alpha_0}+2\frac{\alpha_1}{\alpha_0}|z|^2\big\Vert^2_{L^2}\\
&& + \left(1-4\frac{\alpha_2\alpha_1}{\alpha_0^2}+8\frac{\alpha_1^2}{\alpha_0^2}|z|^2 \right)\Vert \Lambda_\omega F_A \big\Vert^2_{L^2}\\
&&+ \delta(E,[\omega],M,{\alpha}),
\end{eqnarray*}
where $ \delta(E,[\omega],M,\alpha)$ is a constant dependant only on topological
constants of $(E,[\omega],M)$ and the triple 
$\alpha=(\alpha_0,\alpha_1,\alpha_2)$. Remark that we can write
\begin{eqnarray*}\Vert \Lambda_\omega F_A \big\Vert^2_{L^2}&=& \Vert \Lambda_\omega F_A-z \big\Vert^2_{L^2}+2z\int_M F_A\wedge \omega^{n-1} - |z|^2Vol_M([\omega]),\\
&=&\Vert \Lambda_\omega F_A-z \big\Vert^2_{L^2} + \delta'(E,M,[\omega]),
\end{eqnarray*}
where $\delta'(E,M,[\omega])$ depends only on $[\omega]$ and the topology of $(E,M)$.
Therefore, we obtain
\begin{eqnarray*}
CYM(g,A)&= &\big\Vert  Scal(g)-2\frac{\alpha_1}{\alpha_0} |F_A|^2- \frac{\alpha_2}{\alpha_0}+2\frac{\alpha_1}{\alpha_0}|z|^2\big\Vert^2_{L^2}\\
&&+  \left(1-4\frac{\alpha_2\alpha_1}{\alpha_0^2}+8\frac{\alpha_1^2}{\alpha_0^2}|z|^2 \right)\Vert \Lambda_\omega F_A-z \big\Vert^2_{L^2} \\
&& +  \delta''(E,[\omega],M,{\alpha}),
\end{eqnarray*}
with $ \delta''(E,[\omega],M,\alpha)$ a topological constant.  We claim that for several choices of $\omega,L,E$ in Theorem \ref{thm}, there exists a solution to the coupled equations (\ref{cpled}) that minimize the $CYM$ functional. \\
Actually, we remark that if one has the inequality \begin{equation}
 \left(1-4\frac{\alpha_2\alpha_1}{\alpha_0^2}+8\frac{\alpha_1^2}{\alpha_0^2}|z|^2 \right)>0,
 \label{ineqalpha}
                                                                                      \end{equation}
 then we get 
$$CYM(g,A)\geq   \delta''(E,[\omega],M,{\alpha}),$$ and the equality is achieved precisely for a solution to (\ref{cpled}).
From (\ref{ratio2}) and (\ref{fix}), we remark by a direct computation that the limit when $k'\rightarrow +\infty$ of the LHS of (\ref{ineqalpha})  is
$$1+ \frac{3}{k_2^2} + \frac{a^2}{k_2^4}  $$
Thus, if we choose  $k'$ large enough (hence $x>0$ small enough), we get the required inequality (\ref{ineqalpha}) and independently of the choice of the other parameters . \\
Another possible  choice is to take the limit $k_2\rightarrow \pm\infty$ and in that case the LHS of (\ref{ineqalpha}) tends to 1.  
We can also do the following choice : $\alpha_2=0$ with $k_1=k_2=a$, $k'=k$, $b=3a$, $x=\frac{1}{2}$ and $s_\Sigma=-2$ which implies again that the solutions of the coupled equations minimize the $CYM$ functional. \\
From our discussion we obtain the following corollary.
\begin{corollary}
 In Theorem \ref{thm}, there exists for  $k'$ large enough or $|k_2|$ large enough  (or for the choice $s_\Sigma=-2$, $b=3a$, $k'=k$) a solution to the coupled equations (\ref{cpled}) that is the absolute minimum of the $CYM$ functional. 
\end{corollary}
Let us discuss now briefly the uniqueness of the solutions we found. We know that at the level of Chern classes,
$$c_1(M)=2H-p^*c_1(\mathcal{O}\oplus \mathcal{L})+p^* c_1(\Sigma),$$
where $H\in |\mathcal{O}_M(1)|$ and $p:M\rightarrow \Sigma$ is the canonical projection of the ruled manifold $M$ to the surface $\Sigma$.
Thus for our last choice above ($\alpha_2=0$ and $\Sigma$ has genus $h=1+k$, $k\in \mathbb{Z}_+^*$), we are under the conditions of \cite[Proposition 3.5.3 (2)]{FG} since $c_1(M)\leq 0$ for $k$ large enough. Thus we have constructed a family $((M,J),E, 2\pi c_1(L))$ with solutions $(g,A)$ to the coupled equations  (\ref{cpled}) such that the associated K\"ahler form $\omega$ is unique in the K\"ahler class $2\pi c_1(L)$.  
\section{Examples of solutions to coupled equations on the total 
space of a projective bundle over a product of 
two Riemann surfaces}
In order to obtain a form $\Omega$ defined by (\ref{symplform}) which is K\"ahler
we would like to construct an example of solutions to (\ref{cpled}) 
where $\alpha_{1}/\alpha_{0} > 0$.
To that end, we change the setting a little bit to gain more 
flexibility.
\medskip

Let us consider a ruled manifold of the form
$M = {\mathbb P}({\mathcal O} \oplus {\mathcal L}) \rightarrow 
\Sigma_{1} \times \Sigma_{2}$,
where $\Sigma_{i}$, $i=1,2$ is a compact Riemann surface, ${\mathcal 
L} = {\mathcal L}_{1} \otimes {\mathcal L}_{2}$, where ${\mathcal L}_{i}$ is a  
holomorphic line bundle of degree $k_{i}\in {\mathbb Z}^{*}$ on $\Sigma_{i}$, 
and 
${\mathcal O}$ is the trivial holomorphic line bundle.
Let $\pm g_i$ be the
K\"ahler metric
on $\Sigma_{i}$ of constant scalar curvature $\pm2s_i$, with K\"ahler form
$\pm\omega_i$, such that
$c_{1}({\mathcal L_{i}}) = [\frac{\omega_{i}}{2 \pi}]$.
If we
denote the genus of $\Sigma_{i}$ by $h_{i}$,
we have the relation $s_i= 2(1-h_{i})/k_{i}$.
Similarly to the previous section, the zero and infinity sections of 
$M \rightarrow \Sigma_{1} \times 
\Sigma_{2}$ are denoted by $E_{0}$ and $E_{\infty}$. Further, $M_{0}$ and $\mathfrak{z}: 
M_{0} \rightarrow (-1,1)$ are defined as before, while $\theta$ now 
satisfies $d\theta = \omega_{1}+ \omega_{2}$.  Let $x_{1} \neq x_{2}$ 
be real numbers such that $0 < |x_{i}| < 1$ and 
$\frac{\omega_{i}}{x_{i}}$ is positive. 
If, again, $\Theta$ is a smooth real function with 
domain containing
$(-1,1)$ and satisfying \eqref{positivity} we now have an admissible 
metric on $M$ which on $M_{0}$ is given by
\begin{equation}\label{metricproduct}
g  =  \frac{1+x_{1} \mathfrak{z}}{x_{1}} g_1 +\frac{1+x_{2} \mathfrak{z}}{x_{2}} g_2
+\frac {d\mathfrak{z}^2}
{\Theta (\mathfrak{z})}+\Theta (\mathfrak{z})\theta^2\,
\end{equation}
with K\"ahler form 
\begin{equation}
\omega =  \frac{1+x_{1} \mathfrak{z}}{x_{1}} \omega_1 +\frac{1+x_{2} \mathfrak{z}}{x_{2}} 
\omega_2
+d\mathfrak{z}\wedge \theta\,, 
\end{equation}
and the complex structure given as in the previous section. If we set
$$\Theta(\mathfrak{z}) = \frac{F(\mathfrak{z})}{(1+x_{1}\mathfrak{z})(1+x_{2}\mathfrak{z})},$$ then the boundary conditions \eqref{positivity} 
now become equivalent to
\begin{eqnarray}\label{boundaryproduct}
&(i)&\ F(\mathfrak{z}) > 0, \quad -1 < \mathfrak{z} <1,\quad \nonumber\\
&(ii)&\ F(\pm 1) = 0,\quad \nonumber  \\
&(iii)&\ F'(\pm 1) = \mp 2(1 \pm x_{1})(1\pm x_{2}).  \label{positivityFproduct}
\end{eqnarray}
and the scalar curvature of $g$ equals
\begin{equation}\label{scalarproduct}
 Scal_{g} = \frac{2 s_{1} x_{1}}{1+x_{1}\mathfrak{z}} + \frac{2 s_{2} 
x_{2}}{1+x_{2}\mathfrak{z}} - \frac{F''(\mathfrak{z})}{(1+x_{1}\mathfrak{z})(1+x_{2}\mathfrak{z})}.
\end{equation}

Consider the $2$-form $\eta = d(\mathfrak{z}\theta)$ on $M_{0}$. By the 
discussion of Section 1.3 in \cite{acgt} this form is well-defined 
and closed on $M$. In fact, $[\eta]$ is the Poincar\'e dual of 
$2\pi[E_{0}+E_{\infty}]$ and $H^{2}(M,{\mathbb R})$ is generated by $[\eta]$ 
and pullbacks from $\Sigma_{1}\times \Sigma_{2}$. Further $[\eta/(4\pi)]$ 
generates $H^{2}(p^{-1}(x), {\mathbb Z})$ where $p^{-1}(x)$ denotes a 
fibre of $p: M \rightarrow 
\Sigma_{1}\times \Sigma_{2}$ while $ [\omega_{i}/(2\pi k_{i})]$ 
(appropriately lifted) is a 
primitive integer class as well.
We may write the K\"ahler class of our admissible metrics as
\begin{equation}
 [\omega] = \frac{2\pi}{x_{1}}[\omega_{1}/(2\pi)] + 
\frac{2\pi}{x_{2}}[\omega_{2}/(2\pi)] + 4 \pi [\eta/(4\pi)].\label{class}
\end{equation}
Notice that as long as $x_1$ and $x_2$ are rational numbers, $ [\frac{\omega}{2\pi}]$ is a rational class and hence by rescaling the K\"ahler metric as necessary, we obtain an integer class and a corresponding line bundle $L$. Such a rescaling would rescale  $\alpha_{1}/\alpha_{0}$ from the coupled equations by the same factor, but would not change its sign and would not change the qualitative properties of e.g. the CYM functional. In what follows we shall therefore ignore this rescaling factor.

Let
$$ \alpha =  (\omega_{1}+\omega_{2})f+ df \wedge \theta,$$
where $f = \frac{1}{(1+x_{1} \mathfrak{z})(1+x_{2}\mathfrak{z})}$.
It is easy to see that $\alpha$ is closed. If  $\langle.,. \rangle$ denotes 
the inner product on $2$-forms induced by the metric $g$, we observe that $\langle \alpha, \omega \rangle =0$ using the following facts
\begin{eqnarray*}
 \langle \frac{1+x_{i}\mathfrak{z}}{x_{i}} \omega_{i}, \frac{1+x_{j}\mathfrak{z}}{x_{j}} \omega_{j} \rangle &=& \delta_{ij},\\
 \langle \frac{1+x_{i}\mathfrak{z}}{x_{i}} \omega_{i}, d\mathfrak{z} \wedge \theta \rangle &=&0,\\
 \langle d\mathfrak{z} \wedge  \theta, d\mathfrak{z} \wedge \theta \rangle &=&1.
\end{eqnarray*}
Now 
$[\alpha] = m_{1} [\omega_{1}/(2\pi)] + 
m_{2} [\omega_{2}/(2\pi)] + n[\eta/(4\pi)]$ for some $n, m_{1}, 
m_{2} \in {\mathbb R}$. By integrating along 
the fibre $C$ of the ruling and along $\Sigma_{i}$ as embedded in 
$E_{0}= \mathfrak{z}^{-1}(1)$, 
we can determine the value of $n$ and $m_{i}$, $i=1,2$. Indeed, if
$[\alpha] = m_{1} [\omega_{1}/(2\pi)] + 
m_{2} [\omega_{2}/(2\pi)] + n[\eta/(4\pi)]$ then (allowing for a 
slight abuse of notation)
$$ \int_{C}[\alpha] = n$$
while
$$\int_{\Sigma_{i}\subset E_{0}} [\alpha] = 
(m_{i}+\frac{n}{2})\int_{\Sigma_{i}}\frac{\omega_{i}}{2\pi} = 
(m_{i}+\frac{n}{2})k_{i}$$
(using that $\eta = \mathfrak{z}d\theta+d\mathfrak{z}\wedge\theta = 
\mathfrak{z}(\omega_{1}+ \omega_{2})+d\mathfrak{z}\wedge\theta$).
On the other hand, since 
$ \alpha =  (\omega_{1}+\omega_{2})f+ df \wedge \theta$, we also have 
that
$$ \int_{C}[\alpha] =\int_{C} df \wedge \theta = 2\pi (f(1)-f(-1))$$
and 
$$\int_{\Sigma_{i}\subset E_{0}} [\alpha] = f(1) 
\int_{\Sigma_{i}}\omega_{i} = 2\pi f(1)k_{i},$$
and so, using that $f = \frac{1}{(1+x_{1} \mathfrak{z})(1+x_{2}\mathfrak{z})}$, we
get 
$$n= \frac{-4\pi (x_{1}+x_{2})}{(1-x_{1}^{2})(1-x_{2}^{2})}$$
while
$$ m_{1}=m_{2} =  \frac{2\pi(1+x_{1}x_{2})}{(1-x_{1}^{2})(1-x_{2}^{2})}.$$

Similarly to the previous section we now define
$\gamma_{a,b} = a \omega +b \alpha$. From the discussion above we see 
that 
\[
\begin{array}{lcl}
[\gamma_{a,b}]& =& 2 \pi \left(\frac{a}{x_{1}}
+ \frac{b(1+x_{1}x_{2})}{(1-x_{1}^{2})(1-x_{2}^{2})}\right) 
[\omega_{1}/(2\pi)]\\
\\
&+& 2 \pi \left(\frac{a}{x_{2}}
+ \frac{b(1+x_{1}x_{2})}{(1-x_{1}^{2})(1-x_{2}^{2})}\right) 
[\omega_{2}/(2\pi)] \\
\\
&+& 4\pi \left(a- \frac{b 
(x_{1}+x_{2})}{(1-x_{1}^{2})(1-x_{2}^{2})}\right) [\eta/(4\pi)].
\end{array}
\]
Given values of $x_{1}$ and $x_{2}$, it is now clear that we can 
choose $a,b \in {\mathbb R}$ such that
$[\frac{\gamma_{a,b}}{2\pi}]$ is an integer class.

Using that for any $(1,1)$ forms $\beta,\delta$, we have
$$\Lambda^{2}_\omega(\beta \wedge \delta) = 2(\langle \beta,\omega 
\rangle \langle \delta, \omega \rangle - \langle \beta,\delta 
\rangle),$$
we calculate that
\begin{equation*}\begin{array}{lll}\Lambda^{2}_\omega(\gamma_{a,b} \wedge\gamma_{a,b}) &=& 12 a^{2}- 
2b^{2}\frac{x_{1}^{2}(1+x_{2}\mathfrak{z})^{2}+ x_{2}^{2}(1+x_{1}\mathfrak{z})^{2} + \left(x_{1}(1+x_{2}\mathfrak{z}) + 
x_{2}(1+x_{1}\mathfrak{z})\right)^{2}}{(1+x_{1}\mathfrak{z})^{4}(1+x_{2}\mathfrak{z})^{4}}. \\
\end{array}\end{equation*}
Since
$\Lambda_\omega \gamma_{a,b} = 3a$, and the scalar curvature is given by 
\eqref{scalarproduct}
the coupled equations \eqref{cpled} are satisfied, in the case 
$\alpha_{0} \neq 0$, if and only if
\[
    \begin{array}{lcl}
F''(\mathfrak{z}) & = & 2 s_{1} x_{1}(1+x_{2}\mathfrak{z}) + 2 s_{2}x_{2}(1+x_{1}\mathfrak{z})\\
\\
&&+
(12a^{2}\frac{\alpha_{1}}{\alpha_{0}} - 
\frac{\alpha_{2}}{\alpha_{0}})(1+x_{1}\mathfrak{z})(1+x_{2}\mathfrak{z}) \\
\\
&&-2 b^{2}\frac{\alpha_1}{\alpha_0}\frac{x_{1}^{2}(1+x_{2}\mathfrak{z})^{2}+ x_{2}^{2}(1+x_{1}\mathfrak{z})^{2} + (x_{1}(1+x_{2}\mathfrak{z}) + 
x_{2}(1+x_{1}\mathfrak{z}))^{2}}{(1+x_{1}\mathfrak{z})^{3}(1+x_{2}\mathfrak{z})^{3}}.
\end{array}
\]
Set
$\kappa_1= (12a^{2}\frac{\alpha_{1}}{\alpha_{0}} - 
\frac{\alpha_{2}}{\alpha_{0}})$ and $\kappa_2= 4 
b^{2}\frac{\alpha_{1}}{\alpha_{0}}$. The the above equation can be 
written as
\begin{equation}
   \begin{array}{lcl}
       F''(\mathfrak{z}) & = & 2 s_{1} x_{1}(1+x_{2}\mathfrak{z}) + 2 s_{2}x_{2}(1+x_{1}\mathfrak{z}) 
       \\
       \\
&&+ \kappa_1(1+x_{1}\mathfrak{z})(1+x_{2}\mathfrak{z}) \\
\\
& &- \kappa_2 
(x_{1}^{3}(1+x_{1}\mathfrak{z})^{-3}-x_{2}^{3}(1+x_{2}\mathfrak{z})^{-3})/(x_{1}-x_{2}).
\end{array}
\end{equation}

Let us fix
\[
\begin{array}{lcl}
      P(t) & = & \int_{-1}^{t}(2 s_{1} x_{1}(1+x_{2}\mathfrak{z}) + 2 
      s_{2}x_{2}(1+x_{1}\mathfrak{z}) )\, d\mathfrak{z}
       \\
&&+ \kappa_1\int_{-1}^{t}(1+x_{1}\mathfrak{z})(1+x_{2}\mathfrak{z}) \, d\mathfrak{z}\\
&&- \kappa_2 \int_{-1}^{t}(x_{1}^{3}
(1+x_{1}\mathfrak{z})^{-3}-x_{2}^{3}(1+x_{2}\mathfrak{z})^{-3})\frac{1}{(x_{1}-x_{2})}\, d\mathfrak{z}\\
&&+ 2(1-x_{1})(1-x_{2}).
\end{array}
\]
Then 
\[ F(\mathfrak{z}) = \int_{-1}^{\mathfrak{z}}P(t)\, dt\] gives  a bona fide solution if and only 
if $\kappa_1$ and $\kappa_2$ are such that
\begin{equation}\label{prodbasePendpoints}
P(1) = -2(1+x_{1})(1+x_{2}),\quad \quad \int_{-1}^{1}P(t)\,dt = 0,
\end{equation}
and $F(\mathfrak{z}) > 0$ for $-1<\mathfrak{z}<1$.
We calculate that
\[
\begin{array}{lcl}
P(t) & = & (2 s_{1} x_{1} + 2 s_{2}x_{2})(t+1) + 
(s_{1}+s_{2})x_{1}x_{2}(t^{2}-1) \\
&&+ 2(1-x_{1})(1-x_{2})\\
&&+ \kappa_1\left(t+1 + (x_{1}+x_{2})\frac{(t^{2}-1)}{2} + x_{1}x_{2}\frac{(t^{3}+1)}{3}\right)\\
&&+ \frac{\kappa_2}{2(x_{1}-x_{2})}\left(x_{1}^{2}(
\frac{1}{(1+x_{1}t)^{2}}-\frac{1}{(1-x_{1})^{2}})-x_{2}^{2}(
\frac{1}{(1+x_{2}t)^{2}}-\frac{1}{(1-x_{2})^{2}})\right).
\end{array}
\]
Now the first equation of \eqref{prodbasePendpoints} becomes
\[
\begin{array}{l}
(1+\frac{x_{1}x_{2}}{3}) \kappa_1
+ \frac{x_{1}^{3}x_{2}^{3} + 2x_{1}^{2}x_{2}^{2}-x_{1}^{2}-x_{1}x_{2}-x_{2}^{2}}
{(1-x_{1}^{2})^{2}(1-x_{2}^{2})^{2}} \kappa_2 = -2(1+x_{1}x_{2}) -
2(s_{1}x_{1}+s_{2}x_{2}),
\end{array}\]
while the second is
\[
\begin{array}{r}
(1+\frac{x_{1}x_{2}}{3}-\frac{x_{1}+x_{2}}{3}) \kappa_1 
+ \frac{x_{1}^{3}x_{2}^{3} + 
2x_{1}^{2}x_{2}^{2}-x_{1}^{2}-x_{1}x_{2}-x_{2}^{2} 
+(x_{1}+x_{2})(2x_{1}^{2}x_{2}^{2}-x_{1}^{2}-x_{2}^{2})}
{(1-x_{1}^{2})^{2}(1-x_{2}^{2})^{2}} \kappa_2\\
= -2(1+x_{1}x_{2} - (x_{1}+x_{2})) - 2(s_{1}x_{1}+ s_{2}x_{2}) +
\frac{2}{3} (s_{1}+s_{2})x_{1}x_{2}.
\end{array}
\]
Thus we have a linear system in the variables $(\kappa_1,\kappa_2)$ with coefficients determined only by $s_1, s_2, x_1$, and $x_2$. In particular, they do not depend on $a$ and $b$. We spot right away that if 
$x_{1}=-x_{2}$, but $s_{1} \neq -s_{2}$ this system is inconsistent. 
If $x_{1}=-x_{2}$ with
$s_1=s_2$, the system has an infinite number of solutions with say $\kappa_2$ as the free parameter. In particular the K\"ahler class determined by $x_1 = -x_2$ has a constant scalar curvature K\"ahler metric ($\kappa_2=0$) in this case.
On the other hand, if $x_1 \neq -x_2$, it is elementary  to check that the linear system has a unique solution $(\kappa_1,\kappa_2)$. In particular, in that case 
\[ \kappa_2=
\frac{2 (1 - {x_1}^2)^2 (1 - {x_2}^2)^2 (6( {x_1} + x_2) - 3 (s_1 {x_1}^2  + s_2 {x_2}^2) + (s_1+ s_2) {x_1}^2 {x_2}^2 )}{3 ({x_1} + 
     {x_2}) (-4 {x_1}^2 - {x_1} {x_2} - {x_1}^3 {x_2} - 4 {x_2}^2 + 8 {x_1}^2 {x_2}^2 - 
     {x_1} {x_2}^3 + 3 {x_1}^3 {x_2}^3)},
\]
where $(-4 {x_1}^2 - {x_1} {x_2} - {x_1}^3 {x_2} - 4 {x_2}^2 + 8 {x_1}^2 {x_2}^2 - 
     {x_1} {x_2}^3 + 3 {x_1}^3 {x_2}^3)$ in the above expression is always less than zero for
   $0<|x_i| < 1$.
Now, using that $s_ix_i < 2$ it is not hard to check that for $0< x_1, x_2 < 1$ (${\mathcal L}$ from $M = {\mathbb P}({\mathcal O} \oplus {\mathcal L}) \rightarrow 
\Sigma_{1} \times \Sigma_{2}$ being positive definite), $\kappa_2$,  and hence $\frac{\alpha_{1}}{\alpha_{0}}$, is never positive. Likewise, $\kappa_2$ is never positive for $-1< x_1, x_2 < 0$.  
Moving forwards we shall, without loss of generality, assume that $0<x_1<1$ and $-1<x_2 <0$.
Unfortunately,  (i) of \eqref{boundaryproduct} is hard to check in general and by experimenting with some examples we discovered that in some cases it is simply not satisfied. This is a situation not unlike the extremal K\"ahler metric situation on e.g. ruled surfaces of higher genus. For now, we shall focus on a few token examples, taking us through the various genera that may occur for $\Sigma_1$ and $\Sigma_2$, where the positivity of $F(\mathfrak{z})$ may be verified directly.
\subsection{Examples}
\subsubsection{An example with $s_{1}=-s_{2}=2$, $x_{1} = 1/2$, $-1< 
x_{2}=x <0$\label{ex1}}
Since $0<x_{1}<1$ and $s_{1}>0$ while $-1<x<0$ and $s_{2}<0$ this corresponds to the case where $\Sigma_{1}$ and 
$\Sigma_{2}$ both have zero genus and $M =  {\mathbb P}({\mathcal O} \oplus 
{\mathcal O}(1,-1)) \rightarrow 
{\mathbb{CP}}^{1} \times{\mathbb{CP}}^{1}$.

It is easy to check that the linear system has a unique solution 
unless $$x=-1/2.$$ If $x=-1/2$, the system simplifies to a unique equation
\[
33 \kappa_1 - 16 \kappa_2 = -198
\]
which has many solutions, one of them $\kappa_2=0$, i.e. a cscK metric (in 
fact K\"ahler-Einstein \cite{KS}) as it is 
confirmed by
Theorem 9 in \cite{acgt}. This theorem together with \cite[Theorem 8]{acgt} also tells us that no 
other values of $x$ correspond to K\"ahler classes admitting a cscK 
metric. If $x \neq -1/2$, we calculate that
\[
\kappa_1 = \frac{6(x-2)(3+2x+7x^{2})}{8+5x + 16x^{2}+x^{3}}
\]
and 
\[
\kappa_2 = \frac{-9(1-x)^{2}(1+x)^{2}(1+2x)}{8+5x + 16x^{2}+x^{3}}.
\]
We then observe that for $x<-1/2$, $\kappa_2 >0$ and hence 
$\frac{\alpha_{1}}{\alpha_{0}}>0$. For instance, for $x=-3/4$, 
\begin{equation*}
\begin{array}{llr}
P(t)&=&\frac{1}{3284 (2 + t)^2 (4 - 3 t)^2} \left(-57636 - 396428 t + 431692 t^2 + 369508 t^3 \right. \\
&&\left.- 304629 t^4 - 
 150804 t^5 + 60291 t^6 + 25839 t^7\right)
\end{array}
\end{equation*}
   We observe that $P(t)$ is positive at $t=-1$, negative at $t=1$, 
   and
   changes sign only once in the interval $-1<t<1$. Since 
   $\int_{-1}^{1}P(t)\, dt =0$, it is therefore clear that
   $F(\mathfrak{z}) = \int_{-1}^{\mathfrak{z}}P(t)\, dt >0$ for $-1<\mathfrak{z}<1$.
   
\begin{example}\label{thm2}
On $M =  {\mathbb P}({\mathcal O} \oplus 
{\mathcal O}(1,-1)) \rightarrow 
{\mathbb{CP}}^{1} \times{\mathbb{CP}}^{1}$, with  complex structure $J$,
there exists integral classes  $L$ and $E$, a K\"ahler metric $\omega\in 2\pi c_1(L)$, 
and a connection $A\in \mathcal{A}_J^{1,1}(E)$
such that the triple $(\omega,J,A)$ is a solution to the coupled equations (\ref{cpled})
and for which constants $(\alpha_0,\alpha_1,\alpha_2)$ satisfy 
$\frac{\alpha_1}{\alpha_0}>0$. Further, the K\"ahler class $2\pi c_1(L)$ admits no constant scalar curvature K\"ahler metric.
\end{example}

\begin{remark}
Note that $M$ is a toric bundle on compact homogeneous manifolds, so two torus invariant K\"ahler metrics can be joined by a smooth geodesic, see \cite[Theorem 4]{Gu3}.
We also remark that $M$ is a standard compact almost homogeneous space with two ends \cite[Theorem 12.1]{Gu2}.
One can apply \cite[Theorem 2 \& 3]{Gu1} to deduce the existence  of an extremal metric in each K\"ahler class. Note that from \cite[Section 12]{Gu2}, we know that the
geodesics on $M$ satisfy all stability principles and in particular  are smooth. 
Unfortunately we cannot apply the results of  \cite[Section 3.5]{FG} that hold only in dimension 2 to deduce the uniqueness of the solution to the coupled equations.  Nevertheless, we conjecture that the constructed solution $(g,A)$ is unique.
\end{remark}

\subsubsection{An example with $x_1=1/2$, $x_2=-1/3$, $s_1=2$, $s_2=0$}
We apply a similar method to Subsection \ref{ex1}, fixing  $\Sigma_1$ with genus 0 (and $k_1=1$)
and $\Sigma_2$ with genus 
1 (and $k_2<0$).
By applying condition (19) from \cite{acgt}, it is easy to check that 
the K\"ahler class corresponding to these parameters does not admit a 
constant scalar curvature K\"ahler metric.
In that case $\kappa_2=\frac{256}{327}>0$ so $\frac{\alpha_{1}}{\alpha_{0}}>0$.
Furthermore 
\begin{equation*}
\begin{array}{llr}
P(t)&=& \frac{1}{1962 (3 - t)^2 (2 + t)^2}\left(12636 - 120588 t - 85289 t^2 + 33646 t^3\right.\\
&& \left. + 24982 t^4 - 5012 t^5 - 
 2033 t^6 + 394 t^7\right)
\end{array}
\end{equation*}
As before, we observe that $P(t)$ is positive at $t=-1$, negative at $t=1$,  and
   changes sign only once in the interval $-1<t<1$. Since 
   $\int_{-1}^{1}P(t)\, dt =0$, it is therefore clear that
   $F(\mathfrak{z}) = \int_{-1}^{\mathfrak{z}}P(t)\, dt >0$ for $-1<\mathfrak{z}<1$.

\begin{example}\label{thm3}
Let  $M =  {\mathbb P}({\mathcal O} \oplus 
{\mathcal L}_1 \otimes {\mathcal L}_2) \rightarrow 
{\mathbb{CP}}^{1} \times T^2$, where ${\mathcal L}_1 = {\mathcal O}(1) \rightarrow {\mathbb{CP}}^{1}$ and ${\mathcal L}_2 \rightarrow T^2$ is a negative line bundle on $T^2$, and let $J$ denote the  complex structure.
Then there exists integral classes  $L$ and $E$, a K\"ahler metric $\omega\in 2\pi c_1(L)$, 
and a connection $A\in \mathcal{A}_J^{1,1}(E)$
such that the triple $(\omega,J,A)$ is a solution to the coupled equations (\ref{cpled})
and for which constants $(\alpha_0,\alpha_1,\alpha_2)$ satisfy 
$\frac{\alpha_1}{\alpha_0}>0$. Further, the K\"ahler class $2\pi c_1(L)$ admits no constant scalar curvature K\"ahler metric.
\end{example}

\subsubsection{An example with $x_1=1/2$, $x_2=-1/3$, $s_1=2$, $s_2=2$}
We apply a similar method to Subsection \ref{ex1}, fixing  $\Sigma_1$ with genus 0 (and $k_1=1$)
and $\Sigma_2$ with genus 
$h_2 > 1$ (and $k_2= 1-h_2$).
By applying condition (19) from \cite{acgt}, it is easy to check that 
the K\"ahler class corresponding to these parameters does not admit a 
constant scalar curvature K\"ahler metric.
In that case $\kappa_2=\frac{608}{327}>0$ so $\frac{\alpha_{1}}{\alpha_{0}}>0$.
Furthermore 
\begin{equation*}
\begin{array}{llr}
P(t)&=& \frac{2}{981 (3 - t)^2 (2 + t)^2}\left(3456 - 25860 t - 21568 t^2 + 3239 t^3\right.\\
&& \left. +6188 t^4 - 319 t^5 - 502 t^6 + 50 t^7\right)
\end{array}
\end{equation*}
As before, we observe that $P(t)$ is positive at $t=-1$, negative at $t=1$,  and
   changes sign only once in the interval $-1<t<1$. Since 
   $\int_{-1}^{1}P(t)\, dt =0$, it is therefore clear that
   $F(\mathfrak{z}) = \int_{-1}^{\mathfrak{z}}P(t)\, dt >0$ for $-1<\mathfrak{z}<1$.

\begin{example}\label{thm4}
Let  $M =  {\mathbb P}({\mathcal O} \oplus 
{\mathcal L}_1 \otimes {\mathcal L}_2) \rightarrow 
{\mathbb{CP}}^{1} \times \Sigma$, where $\Sigma$ is a Riemann surface of genus at least 2,
${\mathcal L}_1 = {\mathcal O}(1) \rightarrow {\mathbb{CP}}^{1}$, and ${\mathcal L}_2 \rightarrow \Sigma$ is ${\mathcal K}_\Sigma^{\frac{-1}{2}}$ tensored by a flat line bundle on $\Sigma$.
If $J$ denotes the  complex structure,
then there exists integral classes  $L$ and $E$, a K\"ahler metric $\omega\in 2\pi c_1(L)$, 
and a connection $A\in \mathcal{A}_J^{1,1}(E)$
such that the triple $(\omega,J,A)$ is a solution to the coupled equations (\ref{cpled})
and for which constants $(\alpha_0,\alpha_1,\alpha_2)$ satisfy 
$\frac{\alpha_1}{\alpha_0}>0$. Further, the K\"ahler class $2\pi c_1(L)$ admits no constant scalar curvature K\"ahler metric.
\end{example}

\subsubsection{An example with $x_1=1/2$, $x_2=-2/5$, $s_1=0$, $s_2=2$}
We apply a similar method to Subsection \ref{ex1}, fixing  $\Sigma_1$ with genus 1 (and $k_1>0$)
and $\Sigma_2$ with genus 
$h_2 > 1$ (and $k_2= 1-h_2$).
By applying condition (19) from \cite{acgt}, it is easy to check that 
the K\"ahler class corresponding to these parameters does not admit a 
constant scalar curvature K\"ahler metric.
In that case $\kappa_2=\frac{1029}{1475}>0$ so $\frac{\alpha_{1}}{\alpha_{0}}>0$.
Furthermore 
\begin{equation*}
\begin{array}{llr}
P(t)&=& \frac{1}{5310 (2 + t)^2 (5 - 2 t)^2}\left(57622 - 777868 t - 363069 t^2 + 225660 t^3\right.\\
&& \left. +108333 t^4 - 16656 t^5 - 7852 t^6 - 368 t^7\right)
\end{array}
\end{equation*}
As before, we observe that $P(t)$ is positive at $t=-1$, negative at $t=1$,  and
   changes sign only once in the interval $-1<t<1$. Since 
   $\int_{-1}^{1}P(t)\, dt =0$, it is therefore clear that
   $F(\mathfrak{z}) = \int_{-1}^{\mathfrak{z}}P(t)\, dt >0$ for $-1<\mathfrak{z}<1$.

\begin{example}\label{thm5}
Let  $M =  {\mathbb P}({\mathcal O} \oplus 
{\mathcal L}_1 \otimes {\mathcal L}_2) \rightarrow 
T^2 \times \Sigma$, where $\Sigma$ is a Riemann surface of genus at least 2, ${\mathcal L}_1  \rightarrow T^2$ is a positive holomorphic line bundle, and ${\mathcal L}_2 \rightarrow \Sigma$ is ${\mathcal K}_\Sigma^{\frac{-1}{2}}$ tensored by a flat line bundle on $\Sigma$.
If $J$ denotes the  complex structure,
then there exists integral classes  $L$ and $E$, a K\"ahler metric $\omega\in 2\pi c_1(L)$, 
and a connection $A\in \mathcal{A}_J^{1,1}(E)$
such that the triple $(\omega,J,A)$ is a solution to the coupled equations (\ref{cpled})
and for which constants $(\alpha_0,\alpha_1,\alpha_2)$ satisfy 
$\frac{\alpha_1}{\alpha_0}>0$. Further, the K\"ahler class $2\pi c_1(L)$ admits no constant scalar curvature K\"ahler metric.
\end{example}

\subsubsection{An example with $x_1=1/2$, $x_2=-4/9$, $s_1=-1$, $s_2=2$}
We apply a similar method to Subsection \ref{ex1}, fixing  $\Sigma_1$ with genus $h_1>1$ (and $k_1=2(1-h_1)$)
and $\Sigma_2$ with genus 
$h_2 > 1$ (and $k_2= 1-h_2$).
By applying condition (19) from \cite{acgt}, it is easy to check that 
the K\"ahler class corresponding to these parameters does not admit a 
constant scalar curvature K\"ahler metric.
In that case $\kappa_2=\frac{71825}{348408}>0$ so $\frac{\alpha_{1}}{\alpha_{0}}>0$.
Furthermore 
\begin{equation*}
\begin{array}{llr}
P(t)&=& \frac{1}{348408 (2 + t)^2 (9 - 4 t)^2}\left(6466113 - 159543216 t - 40474082 t^2 \right.\\ 
&& \left.  54232672 t^3 +11937913 t^4 - 1961120 t^5 - 731312 t^6 - 579968 t^7\right)
\end{array}
\end{equation*}
As before, we observe that $P(t)$ is positive at $t=-1$, negative at $t=1$,  and
   changes sign only once in the interval $-1<t<1$. Since 
   $\int_{-1}^{1}P(t)\, dt =0$, it is therefore clear that
   $F(\mathfrak{z}) = \int_{-1}^{\mathfrak{z}}P(t)\, dt >0$ for $-1<\mathfrak{z}<1$.

\begin{example}\label{thm6}
Let  $M =  {\mathbb P}({\mathcal O} \oplus 
{\mathcal L}_1 \otimes {\mathcal L}_2) \rightarrow 
\Sigma_1\times \Sigma_2$, where $\Sigma_i$ is a Riemann surface of genus at least 2, ${\mathcal L}_1 \rightarrow \Sigma_1$ is ${\mathcal K}_{\Sigma_1}$ tensored by a flat line bundle on $\Sigma_1$, and ${\mathcal L}_2 \rightarrow \Sigma_2$ is ${\mathcal K}_{\Sigma_2}^{\frac{-1}{2}}$ tensored by a flat line bundle on $\Sigma_2$.
If $J$ denotes the  complex structure,
then there exists integral classes  $L$ and $E$, a K\"ahler metric $\omega\in 2\pi c_1(L)$, 
and a connection $A\in \mathcal{A}_J^{1,1}(E)$
such that the triple $(\omega,J,A)$ is a solution to the coupled equations (\ref{cpled})
and for which constants $(\alpha_0,\alpha_1,\alpha_2)$ satisfy 
$\frac{\alpha_1}{\alpha_0}>0$. Further, the K\"ahler class $2\pi c_1(L)$ admits no constant scalar curvature K\"ahler metric.
\end{example}

The computations in the beginning of Section \ref{CYMsection} are still valid in the present setting and using that $\frac{\alpha_{1}}{\alpha_{0}} = \frac{\kappa_2}{4b^2}$ and
$\frac{\alpha_{2}}{\alpha_{0}} = \frac{3a^2}{b^2}\kappa_2 - \kappa_1$,
and $z=3a$,
we get that 
\begin{eqnarray*}
CYM(g,A)&= &\big\Vert  Scal(g)-2\frac{\alpha_1}{\alpha_0} |F_A|^2- \frac{\alpha_2}{\alpha_0}+2\frac{\alpha_1}{\alpha_0}|z|^2\big\Vert^2_{L^2}\\
&&+  \left(1+\frac{\kappa_1\kappa_2+\frac{3}{2}\frac{a^2}{b^2}\kappa_2^2}{b^2} \right)\Vert \Lambda_\omega F_A-z \big\Vert^2_{L^2} \\
&& +  \delta''(E,[\omega],M,{\alpha}).
\end{eqnarray*}
Since, as we mentioned earlier, $\kappa_1$ and $\kappa_2$ depend only on $x_1,x_2, s_1$, and $s_2$,  we may conclude that for $|b|$ sufficiently large the solutions we constructed above to the coupled system are actually minima of the CYM functional.

\end{document}